\documentclass[10pt]{article}
\usepackage{setspace}

\usepackage{amsthm,amsmath,amssymb}
\usepackage{multirow}
\usepackage{rotating}

\newtheorem{prp}{Proposition}
\newtheorem{cor}{Corollary}
\newtheorem{lem}{Lemma}
\newtheorem{thm}{Theorem}

\begin{document}

\begin{center}
{\Large \bf The set of autotopisms of partial Latin squares.}

\vspace{0.5cm}

{\large Falc\'on, R. M.}

{\small Department of Applied Mathematics I. \\
School of Building Engineering. University of Seville.
Avda. Reina Mercedes 4 A, 41012 - Seville, Spain.
E-mail: {\em rafalgan@us.es}}
\end{center}

\begin{center} {\large \bf Abstract} \end{center}

Symmetries of a partial Latin square are determined by its autotopism group. Analogously to the case of Latin squares, given an isotopism $\Theta$, the cardinality of the set $\mathcal{PLS}_{\Theta}$ of partial Latin squares which are invariant under $\Theta$ only depends on the conjugacy class of the latter, or, equivalently, on its cycle structure. In the current paper, the cycle structures of the set of autotopisms of partial Latin squares are characterized and several related properties studied. It is also seen that the cycle structure of $\Theta$ determines the possible sizes of the elements of $\mathcal{PLS}_{\Theta}$ and the number of those partial Latin squares of this set with a given size. Finally, it is generalized the traditional notion of partial Latin square completable to a Latin square.

\vspace{0.2cm}

\noindent{\bf MSC 2000:} 05B15, 20N05, 20D45.

\noindent{\bf Keywords:} Partial Latin square, autotopism group, cycle structure.

\vspace{0.3cm}

\section{Introduction.}

Every permutation $\pi$ of the symmetric group $S_n$ can be univocally decomposed into product of disjoint cycles. Let $n_{\pi}$ be the number of these cycles. The numbers $\lambda_i^{\pi}$ of cycles of length $i$ in this decomposition determine its {\em cycle structure} as the expression $z_{\pi}=n^{\lambda_n^{\pi}}\cdot...\cdot 2^{\lambda_2^{\pi}}\cdot 1^{\lambda_1^{\pi}}$, where any term of the form $i^0$ is omitted and any term of the form $i^1$ is replaced by $i$. The cardinality of the set $\mathcal{CS}_n$ of possible cycle structures of $S_n$ is equal to the number $p(n)$ of partitions of $n$. Two permutations are conjugate if and only if they have the same cycle structure. Given $\pi\in S_n$, let $\lambda_{\pi}$ and $\pi_{\infty}$ be respectively its length and the union of its $1$-cycles written in natural order. Hereafter, we will suppose $\pi$ to be represented by following their univocal decomposition into a product $\pi_1\pi_2...\pi_{n_{\pi}}$ of disjoint cycles in order of decreasing length, where each cycle $\pi_i$ is written as $(p_{i,1}p_{i,2}...p_{i,\lambda_{\pi_i}})$, with $p_{i,1}=\min_j\{p_{i,j}\}$ and where $p_{i,1}<p_{j,1}$ whenever $i<j$ and $\lambda_{\pi_i}=\lambda_{\pi_j}$. Finally, given $a\in [n]=\{1,2,...,n\}$, we will write $a\in\pi_i$ if there exists $j\in [\lambda_{\pi_i}]$ such that $a=p_{i,j}$. Analogously, $a\in\pi_{\infty}$ will mean $\pi(a)=a$.

A {\em Latin square} of order $n$ is an $n \times n$ array with elements chosen from a set of $n$ distinct symbols such that each symbol occurs precisely once in each row and each column. Hereafter, $[n]$ will be assumed to be this set of symbols and $\mathcal{LS}_n$ will denote the set of Latin squares of order $n$. Given $L=(l_{rc})\in \mathcal{LS}_n$, its {\em orthogonal representation} $O(L)$ is the set of $n^2$ triples $\{(r,c,l_{rc})\,\mid\, r,c\in [n]\}$ defined by the rows $r$, columns $c$ and symbols $l_{rc}$ of $L$. This set verifies the {\em Latin square condition}, i.e., given two triples of $O(L)$ which coincide in two components, then the third component is also the same. Given $\pi\in S_3$, it is defined the Latin square $L^{\pi}$ such that $O(L^{\pi})=\{(l_{\pi(1)},l_{\pi(2)},l_{\pi(3)})\,\mid\, (l_1,l_2,l_3)\in O(L)\}$, which is said to be {\em parastrophic} to $L$. Permutations of rows, columns and symbols also give rise to new Latin squares. Specifically, given three permutations $\alpha,\beta,\gamma$ of the symmetric group $S_n$, the triple $\Theta=(\alpha,\beta,\gamma)\in \mathfrak{I}_n=S_n^3$ is an {\em isotopism} of Latin squares and $L^{\Theta}$ is said to be {\em isotopic} to $L$, where $O(L^{\Theta})=\{(\alpha(r),\beta(c),\gamma(s))\,\mid\, (r, c, s)\in O(L)\}$. To be isotopic is an equivalence relation, which will be denoted by $\sim$, and the set of Latin squares isotopic to $L$ is its {\em isotopism class} $[L]$. The number of Latin squares and isotopism classes of $\mathcal{LS}_n$ are known for $n\leq 11$ \cite{McKay05, Hulpke11}. A list of representatives of isotopism classes for $n\leq 8$ is given in \cite{McKay}. The {\em cycle structure} of $\Theta$ is the triple $z_{\Theta}=(z_{\alpha},z_{\beta},z_{\gamma})$. Hereafter, given a subset $S\subseteq \mathfrak{I}_n$, $\mathcal{CS}_S$ will denote the set of cycle structures of the elements of $S$. Given $z=(z_1,z_2,z_3)\in \mathcal{CS}_S$, where $z_i=n^{z_{in}}\cdot...\cdot 2^{z_{i2}}\cdot 1^{z_{i1}}$, then $n_{z_i}$ will denote the number of cycles of $z_i$, i.e., $n_{z_i}=\sum_{j\in [n]} z_{ij}$. Finally, the {\em parastrophic class} of $z$ is the set $[z]=\{z^{\pi}=(z_{\pi(1)},z_{\pi(2)},z_{\pi(3)}\,\mid\, \pi\in S_3\}$.

If $L^{\Theta}=L$, then $\Theta$ is said to be an {\em autotopism} of $L$. If $\alpha=\beta=\gamma$, then $\Theta$ is an {\em automorphism} of $L$ and $\Theta=\alpha$ is written instead of $(\alpha,\alpha,\alpha)$. Let $\mathcal{LS}_\Theta, \Delta(\Theta)$, $\mathfrak{A}_n$ and $\mathcal{A}_n$ denote respectively the set of Latin squares which have $\Theta$ as an autotopism, its cardinality and the sets of autotopisms and automorphisms of at least one Latin square of order $n$. Necessary conditions for an isotopism to be an autotopism have been given in \cite{Sade68, McKay05, FalconAC, Stones11} and $\mathcal{A}_n$ has been studied in \cite{Wanless04, Bryant09, Kerby10a, Kerby10b, Stones11}. If $\pi\in S_3$ and $L\in\mathcal{LS}_{\Theta}$, then $L^{\pi}\in \mathcal{LS}_{\Theta^{\pi}}$, so permutations on the components of $\Theta$ preserve $\Delta(\Theta)$. Moreover, this cardinality only depends on the conjugacy class of $\Theta$ \cite{FalconAC} or, equivalently, on its cycle structure, so we will also denote it by $\Delta(z_{\Theta})$. A classification of ${\mathcal{CS}_{\mathfrak{A}}}_n$ is known for all $n\leq 17$ \cite{FalconAC, Stones11}. Given $z\in{\mathcal{CS}_{\mathfrak{A}}}_n$, let $\mathfrak{I}_z=\{\Theta\in\mathfrak{I}_n\,\mid\, z_{\Theta}=z\}$ and $\mathcal{LS}_z=\bigcup_{\Theta\in\mathfrak{I}_z}\mathcal{LS}_{\Theta}$.

An {\em incidence structure} is a triple $(P,B,I)$, where $P$ and $B$ respectively are finite sets of {\em points} and {\em blocks} and $I\subseteq P\times B$ is an {\em incidence relation}. It is {\em r-uniform} if every block contains $r$ points and it is {\em s-regular} if every point is exactly on the same number of blocks. Two blocks are {\em equivalent} if they contain the same set of points and the {\em multiplicity} of a block is the size of its equivalence class. In the study of Latin squares, it can be defined a natural incidence relation $I_n$ between $\mathcal{LS}_n$ and $\mathfrak{I}_n$, where, given $L\in \mathcal{LS}_n$ and $\Theta\in \mathfrak{I}_n$, then $(L,\Theta)\in I_n$ if and only if $L\in \mathcal{LS}_{\Theta}$. So, given $z\in {\mathcal{CS}_{\mathfrak{A}}}_n$, the triple $(\mathcal{LS}_z,\mathfrak{I}_z,I_n)$ is a $\Delta(z)$-uniform incidence structure such that every block have the same multiplicity \cite{FalconHMAMS}. Moreover, given $L\in\mathcal{LS}_n$, the triple $([L],\mathfrak{I}_z,I_n)$ is a uniform and regular incidence structure, where every block contains $\Delta_{[L]}(z)$ elements, whose exact value is known for order up to $6$.

Although a general expression for the values of $\Delta(z)$ and $\Delta_{[L]}(z)$ remains unknown, some general and explicit formulas have been given for the former \cite{Laywine81, Laywine85, Falcon06, Stones10} and Gr\"obner bases have been used to know its exact value for all autotopisms of Latin squares of order up to $7$ \cite{FalconMartinJSC07}. For higher orders, Gr\"obner bases have problem with the exponential growth of data storage and the time of computation, for which the use of new combinatorial tools seems to be the key. So, for example, given $\Theta=(\alpha,\beta,\gamma)\in \mathfrak{I}_z$, Gr\"obner bases were used in \cite{FalconMartinEACA08} to obtain the value of $\Delta(z)$ for the majority of the cycle structures of autotopisms of $\mathfrak{A}_8$ and $\mathfrak{A}_9$, by solving the linear equation system formed after adding the constraints $x_{rcs}=x_{\alpha(r)\beta(c)\gamma(s)}$, for all $r, c, s\in [n]$, to those related to the {\em planar 3-index assignment problem} \cite{Euler86}:
$$\begin{array}{l}\min \sum_{r, c, s\in [n]} w_{rcs}\cdot x_{rcs},\\
\begin{array}{ll}\text{subject to } & \sum_{r\in [n]}x_{rcs}=1, \forall c, s\in [n],\\
         \ & \sum_{c\in [n]}x_{rcs}=1, \forall r, s\in [n],\\
         \ & \sum_{s\in [n]}x_{rcs}=1, \forall r, c\in [n],\\
         \ & x_{rcs}\in \{0,1\}, \forall r, c, s \in [n],\end{array}\end{array} \hspace{2cm} (3PAP_n)$$
where $w_{rcs}$ are real weights for all $r, c, s \in [n]$ and whose set of feasible solutions are in $1-1$ correspondence with $\mathcal{LS}_n$ if we define the Latin square $L=(l_{rc})$ such that $l_{rc}=s$ if and only if $x_{rcs}=1$.

All the previous concepts can be naturally extended to {\em partial Latin squares}, i.e., square arrays with elements chosen from a set of $n$ symbols, such that each symbol occurs at most once in each row and in each column. The {\em size} of a partial Latin square $P$ is the number of its non-blank cells and will be denoted by $|P|$. Let $\mathcal{PLS}_n$ and $\mathcal{PLS}_{n,s}$ denote respectively the set of non-empty partial Latin squares of order $n$ and its subset of arrays of size $s$. An upper bound of the elements of $\mathcal{PLS}_{n,s}$ is given in \cite{Ghandehari05}. The {\em orthogonal representation} of $P\in\mathcal{PLS}_n$ is the set $O(P)$ of $|P|$ triples related to the non-blank cells of $P$. Parastrophic partial Latin squares have therefore the same size. Given $\Theta=(\alpha,\beta,\gamma)\in\mathfrak{I}_n$, it is defined the partial Latin square $P^{\Theta}$ such that $O(P^{\Theta})=\{(\alpha(r),\beta(c),\gamma(s))\,\mid\, (r,c,s)\in O(P)\}$, which is said to be {\em isotopic} to $P$ and $[P]$ will denote its {\em isotopism class}. So, $|P^{\Theta}|=|P|$. $\Theta$ is said to be an {\em autotopism} of $P$ if $P^{\Theta}=P$. Let ${\mathfrak{A}_{\mathcal{P}}}_n$ and $\mathcal{PLS}_{\Theta}$ denote respectively the set of autotopisms of $\mathcal{PLS}_n$ and that of non-empty partial Latin squares which have $\Theta$ as an autotopism. Besides, given $z\in {\mathcal{CS}_{\mathfrak{A}_{\mathcal{P}}}}_n$,  $\mathcal{PLS}_z$ will denote the set $\bigcup_{\Theta\in\mathfrak{I}_z} \mathcal{PLS}_{\Theta}$.

A partial Latin square $P\in\mathcal{PLS}_n$ can be {\em completed} to a Latin square $L\in \mathcal{LS}_n$ if $O(P)\subseteq O(L)$. Given $\Theta\in{\mathcal{CS}_{\mathfrak{A}}}_n$, the subset of $\mathcal{LS}_{\Theta}$ of Latin squares to which $P$ can be completed is denoted by $\mathcal{LS}_{\Theta,P}$. The computation of $\Delta(z)$ can be then simplified \cite{FalconMartinJSC07} if a multiplicative factor $c_P\in\mathbb{N}$ is found such that $\Delta(z)=c_P\cdot |\mathcal{LS}_{\Theta,P}|$. Although this factor, which is called {\em $P$-coefficient of symmetry of $\Theta$}, becomes crucial in the processing of high orders, none exhaustive study has been developed in this regard. Indeed, a comprehensive analysis of ${\mathfrak{A}_{\mathcal{P}}}_n$ and $\mathcal{PLS}_{\Theta}$ has not been properly done until now.

The present paper deals with this last question. It is organized as follows: In Section 2, the set ${\mathcal{CS}_{\mathfrak{A}_{\mathcal{P}}}}_n$ will be characterized and several related results exposed. In Section 3, given $z\in {\mathcal{CS}_{\mathfrak{A}_{\mathcal{P}}}}_n$, it will be dealt with the possible sizes of a partial Latin square $P\in \mathcal{PLS}_z$. In Section 4, analogously to the case of Latin squares, it will be proven that the number of partial Latin squares related to an autotopism only depends on the cycle structure of the latter, in such a way that the elements of $\mathcal{PLS}_n$ and $\mathfrak{I}_n$ can be respectively considered as points and blocks of incidence structures whose uniformity and regularity will be studied. Moreover, new constraints will be imposed to the $3PAP_n$ in order to obtain the set $\mathcal{PLS}_{z,s}$ of partial Latin squares of size $s\in [n]$ related to an autotopism of cycle structure $z\in {\mathcal{CS}_{\mathfrak{A}_{\mathcal{P}}}}_n$. Besides, by using Gr\"obner bases, its cardinality $\Delta_s(z)$ will be obtained for $n\leq 4$. Finally, in Section 5, a theoretical ground for the coefficient of symmetry will be exposed. Specifically, given $\Theta\in\mathfrak{I}_n$, it will be studied the set of partial Latin squares of $\mathcal{PLS}_{\Theta}$ which can be completed to Latin squares of $\mathcal{LS}_{\Theta}$.

\section{The set ${\mathcal{CS}_{\mathfrak{A}_{\mathcal{P}}}}_n$.}

Autotopisms of partial Latin squares are univocally determined by their cycle structures:

\begin{lem} \label{lem0a} $\Theta\in {\mathfrak{A}_{\mathcal{P}}}_n$ if and only if $z_{\Theta}\in {\mathcal{CS}_{\mathfrak{A}_{\mathcal{P}}}}_n$.
\end{lem}

{\bf Proof.} The necessary condition holds by definition of ${\mathcal{CS}_{\mathfrak{A}_{\mathcal{P}}}}_n$. Now, if $z\in{\mathcal{CS}_{\mathfrak{A}_{\mathcal{P}}}}_n$, then there must exist $\Theta_0\in \mathfrak{I}_z$ and $P\in \mathcal{PLS}_{\Theta_0}$. So, given $\Theta\in\mathfrak{I}_z$, $\Theta$ and $\Theta_0$ are conjugate and, therefore, there exists $\Theta'\in\mathfrak{I}_n$ such that $\Theta=\Theta'\Theta_0\Theta'^{-1}$. As a consequence, $P^{\Theta'}\in \mathcal{PLS}_{\Theta}$ and $\Theta\in{\mathfrak{A}_{\mathcal{P}}}_n$. \hfill $\Box$

\vspace{0.5cm}

Let us define the set:
$$\mathrm{LCM}_n=\{(i,j,k)\in [n]^3\,\mid\,\mathrm{lcm}(i,j)=\mathrm{lcm}(i,k)=\mathrm{lcm}(j,k)=\mathrm{lcm}(i,j,k)\}.$$

The next result characterizes the set ${\mathcal{CS}_{\mathfrak{A}_{\mathcal{P}}}}_n$ and can be considered as an immediate generalization for partial Latin squares of the necessary condition given by Stones, Vojt\v{e}chovsk\'y and Wanless in \cite{Stones11} for membership in $\mathfrak{A}_n$:

\begin{lem} \label{lem0} Given $z=(z_1,z_2,z_3)\in \mathcal{CS}_{\mathfrak{I}_n}$, it is $z\in {\mathcal{CS}_{\mathfrak{A}_{\mathcal{P}}}}_n$ if and only if there exists $(i,j,k)\in \mathrm{LCM}_n$ such that $z_{1i}\cdot z_{2j}\cdot z_{3k}>0$.
\end{lem}

{\bf Proof.} If $z\in {\mathcal{CS}_{\mathfrak{A}_{\mathcal{P}}}}_n$, then there must exist $\Theta=(\alpha,\beta,\gamma)\in\mathfrak{I}_z$ and $P\in\mathcal{PLS}_{\Theta}$. Given $(r,c,s)\in O(P)$, let $(u,v,w)\in [n_{\alpha}]\times[n_{\beta}]\times[n_\gamma]$ be such that $r\in \alpha_u$, $c\in\beta_v$ and $s\in\gamma_w$. Since $\Theta$ is an autotopism of $P$, it must be $(\alpha_u^t(r),\beta_v^t(c),\gamma_w^t(s))\in O(P)$, for all $t\in\mathbb{N}$. The necessary condition is then a consequence of the Latin square condition, by considering $i,j,k$ to be, respectively, the lengths of $\alpha_u$, $\beta_v$ and $\gamma_w$.

To prove the converse, let $\Theta=(\alpha,\beta,\gamma)\in\mathfrak{I}_z$ and let $\alpha_u, \beta_v$ and $\gamma_w$ be, respectively, $i$-, $j$- and $k$-cycles of $\alpha$, $\beta$ and $\gamma$. Let $r, c, s$ be, respectively, elements of $\alpha_u, \beta_v$ and $\gamma_w$. The set of triples $\{(\alpha_u^t(r),\beta_v^t(c),\gamma_w^t(s))\,\mid\, t\in [\mathrm{lcm}(i,j,k)]\}$ verifies the Latin square condition because of being $(i,j,k)\in \mathrm{LCM}_n$ and, therefore, it is the orthogonal representation of a partial Latin square $P\in \mathcal{PLS}_{\Theta}$. \hfill $\Box$

\vspace{0.5cm}

Given $n>1$, Lemma \ref{lem0} implies ${\mathfrak{A}_{\mathcal{P}}}_n$ to be a proper subset of $\mathfrak{I}_n$, because, for instance, $(1^2, 1^2, n^1)\not\in {\mathcal{CS}_{\mathfrak{A}_{\mathcal{P}}}}_n$. Analogously, $\mathfrak{A}_n$ is a proper subset of ${\mathfrak{A}_{\mathcal{P}}}_n$, because, for example, $(2,2,2)\in{\mathcal{CS}_{\mathfrak{A}_{\mathcal{P}}}}_2$ and $(2\cdot 1^{n-2},2\cdot 1^{n-2},1^n)\in{\mathcal{CS}_{\mathfrak{A}_{\mathcal{P}}}}_n$, for $n>2$, but neither of them are cycle structures of an autotopism of Latin square. Thus, the next claim is verified:

\begin{prp}\label{prp0} ${{\mathfrak{A}}}_n\subset {{\mathfrak{A}_{\mathcal{P}}}}_n\subset {{\mathfrak{I}}}_n$, $\forall n>1$. \hfill $\Box$
\end{prp}

\vspace{0.5cm}

Moreover, ${{\mathfrak{A}_{\mathcal{P}}}}_n$ and ${{\mathfrak{I}}}_n$ can be identified when $n$ tends to infinity. To see it, it is enough to assure that the cardinalities of the sets of their cycle structures coincide in the limit, which will be proven in Theorem \ref{thm0}. Previously, although ${\mathcal{CS}_{\mathfrak{A}_{\mathcal{P}}}}_n$ can be explicitly obtained for any order $n\in\mathbb{N}$ by implementing Lemma \ref{lem0} in a computer procedure, a lower bound of its cardinality will be determined by studying the following sets which partition $\mathcal{CS}_n$:
$$\mathcal{CS}_{n,m}=\{n^{z_n}\cdot ...\cdot 2^{z_2}\cdot 1^{z_1}\in \mathcal{CS}_n\,\mid\, z_m>0 \text{ and } z_i=0, \forall i\in [m-1]\},$$
where $m\in[n]$. The following results hold:

\begin{lem}\label{lemCS} $|\mathcal{CS}_{n,m}|=\begin{cases} 1, \text{ if } m=n,\\
                                                             0, \text{ if } m\in \{\lceil\frac n2 \rceil,...,n-1\},\\
                                                             p(n-m) - \sum_{i=1}^{m-1} |\mathcal{CS}_{n-m,i}|, \text{ otherwise}.\end{cases}$
\end{lem}

{\bf Proof.} The cases $m\geq \lceil\frac n2 \rceil$ are straightforward verified. Let $m\leq \lfloor\frac n2 \rfloor$. Given $z\in \mathcal{CS}_{n,m}$, it is $z_m>0$ and, therefore, it must be $z_{n-m+i}=0$, for all $i\in [m]$. So, we can define $z'\in \mathcal{CS}_{n-m}$ such that $z'_i=z$, for all $i\in [n-m]\setminus\{m\}$ and $z'_m=z_m-1$. Specifically, since $z\in \mathcal{CS}_{n,m}$, it must be $z'_i=0$, for all $i\in [m-1]$. Thus, $z'\in \bigcup_{i=m}^{n-m}\mathcal{CS}_{n-m,i}=\mathcal{CS}_{n-m}\setminus \bigcup_{i=1}^{m-1}\mathcal{CS}_{n-m,i}$ and the claim is verified. \hfill $\Box$

\vspace{0.5cm}

\begin{prp}\label{prpCS} $|{\mathcal{CS}_{\mathfrak{A}_{\mathcal{P}}}}_n|\geq \sum_{(i,j,k)\in \mathrm{LCM}_n} |\mathcal{CS}_{n,i}|\cdot |\mathcal{CS}_{n,j}|\cdot |\mathcal{CS}_{n,k}|$.
\end{prp}

{\bf Proof.} Since the sets $\mathcal{CS}_{n,m}$ constitute a partition of $\mathcal{CS}_n$, the result is consequence of Lemma \ref{lem0}.\hfill $\Box$

\vspace{0.5cm}

\begin{thm}\label{thm0} $\lim_{n\rightarrow\infty} \frac {|{\mathcal{CS}_{\mathfrak{A}_{\mathcal{P}}}}_n|}{|{\mathcal{CS}_{\mathfrak{I}}}_n|}=1$.
\end{thm}

{\bf Proof.} Since $(1,1,1)\in \mathrm{LCM}_n$, Proposition \ref{prpCS} implies that $|{\mathcal{CS}_{\mathfrak{A}_{\mathcal{P}}}}_n|\geq |\mathcal{CS}_{n,1}|^3$ and, therefore, from Lemma \ref{lemCS}, $|{\mathcal{CS}_{\mathfrak{A}_{\mathcal{P}}}}_n|\geq p(n-1)^3$. Now, since $p(n)$ is equivalent to $\frac {e^{\pi\sqrt{2n/3}}}{4n\sqrt{3}}$ when $n$ tends to infinity \cite{Hardy18} and $|{\mathcal{CS}_{\mathfrak{A}_{\mathcal{P}}}}_n|\leq |{\mathcal{CS}_{\mathfrak{I}}}_n|=p(n)^3$, then:
$$1\geq\lim_{n\rightarrow\infty} \frac {|{\mathcal{CS}_{\mathfrak{A}_{\mathcal{P}}}}_n|}{|{\mathcal{CS}_{\mathfrak{I}}}_n|} \geq \lim_{n\rightarrow\infty} \frac {p(n-1)^3}{p(n)^3}=1.$$ \hfill $\Box$

\vspace{0.5cm}

For $n\leq 17$, Table 1 shows the values $|\mathcal{CS}_{n,m}|$ and $|{\mathcal{CS}_{\mathfrak{A}_{\mathcal{P}}}}_n|$, where $m\leq \lfloor \frac n2\rfloor$, in comparison with those of $|{\mathcal{CS}_{\mathfrak{A}}}_n|$, which can be obtained by using the classification given in \cite{FalconAC, Stones11}. The number $|[{\mathcal{CS}_{\mathfrak{A}_{\mathcal{P}}}}_n]|$ of parastrophic classes of ${\mathcal{CS}_{\mathfrak{A}_{\mathcal{P}}}}_n$ is also shown.

\begin{table}\label{t1}
\begin{center}
\begin{tabular}{|c|c|c|c|c|c|c|c|c|c|c|c|}\hline
\multirow{3}{*}{$n$} & \multirow{3}{*}{$|{\mathcal{CS}_{\mathfrak{A}}}_n|$} & \multicolumn{8}{|c|}{$|\mathcal{CS}_{n,m}|$}& \multirow{3}{*}{$|{\mathcal{CS}_{\mathfrak{A}_{\mathcal{P}}}}_n|$} & \multirow{3}{*}{$|[{\mathcal{CS}_{\mathfrak{A}_{\mathcal{P}}}}_n]|$}\\\cline{3-10}
\ & \ & \multicolumn{8}{|c|}{$m$}& \ & \ \\ \cline{3-10}
\ & & 1 & 2 & 3 & 4 & 5 & 6 & 7 & 8& \ & \\\hline
1 & 1 & \multicolumn{8}{|c|}{} & 1 & 1\\ \cline{1-3}  \cline{11-12}
2 & 4 & 1 & \multicolumn{7}{|c|}{}& 5 & 3\\\cline{1-3}  \cline{11-12}
3 & 6 & 2 & \multicolumn{7}{|c|}{} & 15 & 7\\\cline{1-4} \cline{11-12}
4 & 19 & 3 & 1 & \multicolumn{6}{|c|}{}& 65 & 22\\\cline{1-4} \cline{11-12}
5 & 8 & 5 & 1 & \multicolumn{6}{|c|}{}& 223 & 60\\\cline{1-5} \cline{11-12}
6 & 45 & 7 & 2 & 1 & \multicolumn{5}{|c|}{}& 869 & 197\\\cline{1-5} \cline{11-12}
7 & 12 & 11 & 2 & 1 & \multicolumn{5}{|c|}{}& 2535 & 526\\\cline{1-6} \cline{11-12}
8 & 87 & 15 & 4 & 1 & 1 & \multicolumn{4}{|c|}{} & 7663 & 1492\\\cline{1-6} \cline{11-12}
9 & 43 & 22 & 4 & 2 & 1 & \multicolumn{4}{|c|}{} & 21156 & 3937\\\cline{1-7} \cline{11-12}
10 & 89 & 30 & 7 & 2 & 1 & 1\ & \multicolumn{3}{|c|}{} & 60264 & 10850\\\cline{1-7} \cline{11-12}
11 & 21 & 42 & 8 & 3 & 1 & 1 & \multicolumn{3}{|c|}{} & 150953 & 26628\\\cline{1-8} \cline{11-12}
12 & 407 & 56 & 12 & 4 & 2 & 1 & 1 & \multicolumn{2}{|c|}{} & 385538 & 66984\\\cline{1-8} \cline{11-12}
13 & 27 & 77 & 14 & 5 & 2 & 1 & 1 & \multicolumn{2}{|c|}{} & 915452 & 157398\\\cline{1-9} \cline{11-12}
14 & 141 & 101 & 21 & 6 & 3 & 1 & 1 & 1 & \ & 2193225 & 374127\\\cline{1-9} \cline{11-12}
15 & 150 & 135 & 24 & 9 & 3 & 2 & 1 & 1 & \ & 4928696 & 836154\\\hline
16 & 503 & 176 & 34 & 10 & 5 & 2 & 1 & 1 & 1 & 11209311 & 1893607\\\hline
17 & 40 & 231 & 41 & 13 & 5 & 3 & 1 & 1 & 1 & 24406191 & 4110132\\\hline
\end{tabular}
\caption{Cardinality of the sets of cycle structures, for $n\leq 17$ and $m\leq \lfloor \frac n2\rfloor$.}
\end{center}
\end{table}

\section{The size of a partial Latin square related to an autotopism.}

Given $z=(z_1,z_2,z_3)\in {\mathcal{CS}_{\mathfrak{A}_{\mathcal{P}}}}_n$ and $\Theta=(\alpha,\beta,\gamma)\in\mathfrak{I}_z$, any partial Latin square $P\in\mathcal{PLS}_{\Theta}$ can be decomposed into $n_{z_1}\cdot n_{z_2}$ blocks $P_{ij}$ whose rows and columns are respectively determined by the elements of the cycle $\alpha_i$ of $\alpha$ and the cycle $\beta_j$ of $\beta$, i.e., $O(P_{ij})=\{(r,c,s)\in O(P)\,\mid\, r\in\alpha_i \text{ and } c\in\beta_j\}$. It will be called the {\em $\Theta$-decomposition of $P$}. Specifically, $z$ determines not only the number of these blocks, but also their possible sizes and, consequently, a pair of bounds for the size of $P$. To see it, let us define the set:
$$\mathrm{LCM}_z=\{(i,j)\in [n]^2\,\mid\, \exists k\in[n] \text{ s.t. } (i,j,k)\in \mathrm{LCM}_n \text{ and } z_{1i}\cdot z_{2j}\cdot z_{3k}>0\}.$$

The following results hold:

\begin{lem}\label{lemR1} Given $z\in {\mathcal{CS}_{\mathfrak{A}_{\mathcal{P}}}}_n$, $\Theta\in \mathfrak{I}_z$, $P\in\mathcal{PLS}_{\Theta}$ and an $i\times j$-block $B$ of the $\Theta$-decomposition of $P$, there exists $\omega_B\in [\mathrm{gcd}(i,j)]\cup\{0\}$ such that  $|B|=\omega_B\cdot \mathrm{lcm}(i,j)$. Specifically, $\omega_B=0$ if $(i,j)\not\in \mathrm{LCM}_z$.
\end{lem}

{\bf Proof.} Analogously to the proof of Lemma \ref{lem0}, the Latin square condition implies $|B|=0$, whenever $(i,j)\not\in \mathrm{LCM}_z$. Besides, given $(r,c,s)\in O(B)$, its orbit by the action of $\Theta=(\alpha,\beta,\gamma)$ is the set of triples $(\alpha^t(r),\beta^t(c),$ $\gamma^t(s))\in O(B)$, for all $t\in [\mathrm{lcm}(i,j)]$. So, the Latin square condition implies $|B|$ to be a multiple of $\mathrm{lcm}(i,j)$. Finally, since there are $i\cdot j$ cells in $B$, the multiplicative factor must be at most $\mathrm{gcd}(i,j)$. \hfill $\Box$

\vspace{0.5cm}

\begin{prp}\label{prpR1} Given $z=(z_1,z_2,z_3)\in {\mathcal{CS}_{\mathfrak{A}_{\mathcal{P}}}}_n$ and $P\in\mathcal{PLS}_z$, it is $\mathfrak{l}_z\leq |P|\leq \mathfrak{u}_z$, where:
$$\mathfrak{l}_z=\min_{(i,j)\in \mathrm{LCM}_z}\{\mathrm{lcm}(i,j)\},$$
{\small $$\mathfrak{u}_z=\min\{\sum_{(i,j)\in \mathrm{LCM}_z} z_{1i}\cdot z_{2j} \cdot i\cdot j, \sum_{(i,k)\in \mathrm{LCM}_{z^{(23)}}} z_{1i}\cdot z_{3k} \cdot i\cdot k, \sum_{(k,j)\in \mathrm{LCM}_{z^{(13)}}} z_{2j}\cdot z_{3k} \cdot j\cdot k\}.$$}
\end{prp}

{\bf Proof.} Let $\Theta\in \mathfrak{I}_z$ be such that $P\in\mathcal{PLS}_{\Theta}$ and let $B$ be a block of the $\Theta$-decomposition of $P$ such that $|B|>0$. From Lemma \ref{lemR1}, if $B$ is an $i\times j$-block, where $(i,j)\in\mathrm{LCM}_z$, then $\mathrm{lcm}(i,j)\leq |B|\leq i\cdot j$ and so, $\mathfrak{l}_z\leq |P|\leq \sum_{(i,j)\in \mathrm{LCM}_z} z_{1i}\cdot z_{2j} \cdot i\cdot j$. Since the size of a partial Latin square is invariant by parastrophism and $P^{\pi}\in \mathcal{PLS}_{\Theta^{\pi}}$ for all $\pi\in S_3$, then $\mathfrak{u}_z$ is an upper bound of $|P|$. \hfill $\Box$

\vspace{0.5cm}

From the previous results, it is deduced that the possible sizes of the elements of $\mathcal{PLS}_z$ must be in the set:
$$\mathrm{Sizes}(z)=\left\{\sum_{(i,j)\in \mathrm{LCM}_z} \omega_{ij} \cdot \mathrm{lcm}(i,j)\leq\mathfrak{u}_z\,\mid\, \omega_{ij}\in [z_{1i}\cdot z_{2j}\cdot \mathrm{gcd}(i,j)]\right\}.$$

As an example, let us consider $z=(6, 3\cdot 2\cdot 1, 4 \cdot 2)\in{\mathcal{CS}_{\mathfrak{A}_{\mathcal{P}}}}_6$ and $\Theta=((123456),(123)(45)(6),(1234)(56))\in\mathcal{PLS}_z$. The $\Theta$-decomposition of any partial Latin square $P\in \mathcal{PLS}_{\Theta}$ is then formed by three blocks, $P_{11}, P_{12}$ and $P_{13}$, whose cells are respectively indicated by the symbols $\cdot$, $*$ and $\circ$ in the following diagram:
{\scriptsize $$\left(\begin{array}{cccccc} \cdot & \cdot & \cdot & * &* &\circ\\
                                      \cdot & \cdot & \cdot & * &* &\circ\\
                                      \cdot & \cdot & \cdot &* &* &\circ\\
                                      \cdot & \cdot & \cdot &* &* &\circ\\
                                      \cdot & \cdot & \cdot &* &* &\circ\\
                                      \cdot & \cdot & \cdot &* &* &\circ\\
\end{array}\right).$$}

Besides, $\mathrm{LCM}_z=\{(6,3)\}$, $\mathrm{LCM}_{z^{(23)}}=\{(6,2)\}$ and $\mathrm{LCM}_{z^{(13)}}=\{(2,3)\}$. So, from Proposition \ref{prpR1}, it must be $6\leq |P|\leq \min\{18,12,6\}=6$. Thus, $\mathrm{Sizes}(z)=\{6\}$ and $|P|=6$. Specifically, there are six possibilities for $P$:
{\scriptsize $$\left(\begin{array}{cccccc} 5 & \cdot & \cdot &\cdot &\cdot &\cdot\\
                              \cdot & 6 & \cdot &\cdot &\cdot &\cdot\\
                              \cdot & \cdot & 5 &\cdot &\cdot &\cdot\\
                              6 & \cdot & \cdot &\cdot &\cdot &\cdot\\
                              \cdot & 5 & \cdot &\cdot &\cdot &\cdot\\
                              \cdot & \cdot & 6 &\cdot &\cdot &\cdot\\
\end{array}\right), \hspace{1cm} \left(\begin{array}{cccccc} 6 & \cdot & \cdot &\cdot &\cdot &\cdot\\
                              \cdot & 5 & \cdot &\cdot &\cdot &\cdot\\
                              \cdot & \cdot & 6 &\cdot &\cdot &\cdot\\
                              5 & \cdot & \cdot &\cdot &\cdot &\cdot\\
                              \cdot & 6 & \cdot &\cdot &\cdot &\cdot\\
                              \cdot & \cdot & 5 &\cdot &\cdot &\cdot\\
\end{array}\right), \hspace{1cm} \left(\begin{array}{cccccc} \cdot & 5 & \cdot &\cdot &\cdot &\cdot\\
                              \cdot & \cdot & 6 &\cdot &\cdot &\cdot\\
                              5 & \cdot & \cdot &\cdot &\cdot &\cdot\\
                              \cdot & 6 & \cdot &\cdot &\cdot &\cdot\\
                              \cdot & \cdot & 5 &\cdot &\cdot &\cdot\\
                              6 & \cdot & \cdot &\cdot &\cdot &\cdot\\
\end{array}\right),$$
$$\left(\begin{array}{cccccc} \cdot & 6 & \cdot &\cdot &\cdot &\cdot\\
                              \cdot & \cdot & 5 &\cdot &\cdot &\cdot\\
                              6 & \cdot & \cdot &\cdot &\cdot &\cdot\\
                              \cdot & 5 & \cdot &\cdot &\cdot &\cdot\\
                              \cdot & \cdot & 6 &\cdot &\cdot &\cdot\\
                              5 & \cdot & \cdot &\cdot &\cdot &\cdot\\
\end{array}\right),\hspace{1cm} \left(\begin{array}{cccccc} \cdot & \cdot & 5 &\cdot &\cdot &\cdot\\
                              6 & \cdot & \cdot &\cdot &\cdot &\cdot\\
                              \cdot & 5 & \cdot &\cdot &\cdot &\cdot\\
                              \cdot & \cdot & 6 &\cdot &\cdot &\cdot\\
                              5 & \cdot & \cdot &\cdot &\cdot &\cdot\\
                              \cdot & 6 & \cdot &\cdot &\cdot &\cdot\\
\end{array}\right),\hspace{1cm} \left(\begin{array}{cccccc} \cdot & \cdot & 6 &\cdot &\cdot &\cdot\\
                              5 & \cdot & \cdot &\cdot &\cdot &\cdot\\
                              \cdot & 6 & \cdot &\cdot &\cdot &\cdot\\
                              \cdot & \cdot & 5 &\cdot &\cdot &\cdot\\
                              6 & \cdot & \cdot &\cdot &\cdot &\cdot\\
                              \cdot & 5 & \cdot &\cdot &\cdot &\cdot\\
\end{array}\right).$$}

\section{The number of partial Latin squares related to an autotopism.}

Given $z\in {\mathcal{CS}_{\mathfrak{A}_{\mathcal{P}}}}_n$, $\Theta\in\mathfrak{I}_z$, $P\in\mathcal{PLS}_{\Theta}$ and $s\in [n]$, let us define the sets:
$$\mathcal{PLS}_{\Theta,[P]}=\mathcal{PLS}_{\Theta}\cap [P], \hspace{1.5cm} \mathcal{PLS}_{\Theta,s}=\mathcal{PLS}_{\Theta}\cap \mathcal{PLS}_{n,s}.$$

In the current section, given $\Theta\in\mathfrak{I}_n$, the cardinality of the set $\mathcal{PLS}_{\Theta}$ will be studied. The following result implies that it only depends on the cycle structure of $\Theta$:

\begin{lem} \label{lemD1} The number of isotopic partial Latin squares related to an autotopism only depends on the parastrophic class of the cycle structure of the latter.
\end{lem}

{\bf Proof.} Let $\Theta_1, \Theta_2\in \mathfrak{I}_n$ and $\pi\in S_3$ be such that $z_{\Theta_1}=z_{\Theta_2^{\pi}}$ and let  $P\in\mathcal{PLS}_n$. Since $\Theta_1$ and $\Theta_2^{\pi}$ are conjugate, there exists $\Theta\in \mathfrak{I}_n$ such that $\Theta_2^{\pi}=\Theta\Theta_1\Theta^{-1}$. Now, given $Q\in\mathcal{PLS}_{\Theta_1,[P]}$, it is $Q^{\Theta}\in \mathcal{PLS}_{\Theta_2^{\pi},[P]}$ and, therefore, $\left(Q^{\Theta}\right)^{\pi^{-1}}\in \mathcal{PLS}_{\Theta_2,[P]}$. So, $|\mathcal{PLS}_{\Theta_1,[P]}|\leq|\mathcal{PLS}_{\Theta_2,[P]}|$. The opposite inequality is similarly proven. \hfill $\Box$

\vspace{0.5cm}

Since the size of a partial Latin square is preserved by isotopism, Lemma \ref{lemD1} implies the following cardinalities to be well-defined:
$$\Delta_{[P]}(z)=\left|\mathcal{PLS}_{\Theta, [P]}\right|,$$
$$\Delta_s(z)=\left|\mathcal{PLS}_{\Theta,s}\right|=\sum_{\scriptsize \begin{array}{c}[Q]\in\mathcal{PLS}_{\Theta}/\sim\\s.t.\ |Q|=s\end{array}}\Delta_{[Q]}(z),$$
$$\Delta_{\mathcal{P}}(z)=\left| \mathcal{PLS}_{\Theta}\right|=\sum_{\scriptsize [Q]\in\mathcal{PLS}_{\Theta}/\sim}\Delta_{[Q]}(z)=\sum_{s\in \mathrm{Sizes}(z)}\Delta_s(z).$$

\vspace{0.5cm}

It can be defined a natural incidence relation $I_{\mathcal{P}_n}$ between $\mathcal{PLS}_n$ and $\mathfrak{I}_n$, where, given $P\in \mathcal{PLS}_n$ and $\Theta\in \mathfrak{I}_n$, then $(P,\Theta)\in I_n$ if and only if $P\in \mathcal{PLS}_{\Theta}$. Besides, let us denote by $\mathfrak{A}_P$ the set of autotopisms of $P$. The following results are then proven:

\begin{prp}\label{prpD1} Let $P\in \mathcal{PLS}_n$ and $z\in {\mathcal{CS}_{\mathfrak{A}_{\mathcal{P}}}}_n$. The triples  $([P],\mathfrak{I}_z,I_{\mathcal{P}_n})$, $(\mathcal{PLS}_{n,s},\mathfrak{I}_z,I_{\mathcal{P}_n})$ and $(\mathcal{PLS}_n,\mathfrak{I}_z,I_{\mathcal{P}_n})$ are, respectively, $\Delta_{[P]}(z)$-, $\Delta_s(z)$- and $\Delta_{\mathcal{P}}(z)$-uniform incidence structures and all its blocks have the same multiplicity. Moreover, the former incidence structure is regular.
\end{prp}

{\bf Proof.} From Lemma \ref{lemD1}, it is enough to study the uniformity and multiplicity of $([P],\mathfrak{I}_z,I_{\mathcal{P}_n})$. Indeed, the uniformity is an immediate consequence of that lemma. Now, in order to see that all the blocks have the same multiplicity, let $\Theta_1,\Theta'_1\in\mathfrak{I}_z$ be such that $\mathcal{PLS}_{\Theta_1,[P]}=\mathcal{PLS}_{\Theta'_1,[P]}$ and let us consider $\Theta_2\in\mathfrak{I}_z$. Let $\Theta,\Theta'\in\mathfrak{I}_n$ be such that $\Theta_1=\Theta\Theta_2\Theta^{-1}$ and $\Theta'_1=\Theta'\Theta_1\Theta'^{-1}$. Then, $\mathcal{PLS}_{\Theta_2,[P]}=\mathcal{PLS}_{\Theta^{-1}\Theta'\Theta,[P]}$, because $Q\in\mathcal{PLS}_{\Theta_2,[P]}\Leftrightarrow Q^{\Theta}\in\mathcal{PLS}_{\Theta_1,[P]}\Leftrightarrow Q^{\Theta'\Theta}\in\mathcal{PLS}_{\Theta'_1,[P]}=\mathcal{PLS}_{\Theta_1,[P]}\Leftrightarrow Q^{\Theta^{-1}\Theta'\Theta}\in\mathcal{PLS}_{\Theta_2,[P]}$. Moreover, $\Theta^{-1}\Theta'\Theta=\Theta_2 \Leftrightarrow \Theta'=\Theta^{-1}\Theta_2\Theta=\Theta_1 \Leftrightarrow \Theta'_1=\Theta_1$. Thus, the arbitrariness of $\Theta_1,\Theta'_1$ and $\Theta_2$ implies the claim about the multiplicity.

Finally, in order to see that $([P],\mathfrak{I}_z,I_{\mathcal{P}_n})$ is a regular incidence structure, let us consider $Q_1, Q_2\in [P]$ and let $\Theta\in\mathfrak{A}_{Q_1}\cap \mathfrak{I}_z$. Since $Q_1$ and $Q_2$ are isotopic, there must exist $\Theta'\in \mathfrak{I}_n$ such that $Q_1^{\Theta'}=Q_2$. So, $\Theta'\Theta\Theta'^{-1}\in \mathfrak{A}_{Q_2}\cap\mathfrak{I}_z$ and, therefore, $|\mathfrak{A}_{Q_1}\cap\mathfrak{I}_z|\leq |\mathfrak{A}_{Q_2}\cap\mathfrak{I}_z|$. The regularity holds because the opposite inequality is analogously proven. \hfill $\Box$

\vspace{0.5cm}

\begin{thm}\label{thmD2} Let $P\in\mathcal{PLS}_n$. If $Q\in [P]$, then $|\mathfrak{A}_Q|=|\mathfrak{A}_P|$ and it coincides with the cardinality of the set $\mathfrak{I}_{P,Q}$ of isotopisms from $P$ to $Q$. Moreover, given $\Theta\in\mathfrak{A}_P$, it is $\mathfrak{A}_Q=\{\Theta'\Theta{\Theta'}^{-1}\,\mid\, \Theta'\in \mathfrak{I}_{P,Q}\}$.
\end{thm}

{\bf Proof.} From Proposition \ref{prpD1}, it is verified that $|\mathfrak{A}_Q|=\sum_{z\in {\mathcal{CS}_{\mathfrak{A}_{\mathcal{P}}}}_n}|\mathfrak{A}_Q\cap\mathfrak{I}_z|=\sum_{z\in {\mathcal{CS}_{\mathfrak{A}_{\mathcal{P}}}}_n}|\mathfrak{A}_P\cap\mathfrak{I}_z|=|\mathfrak{A}_P|$. Besides, given $\Theta\in\mathfrak{A}_P$, it is $\Theta'\Theta{\Theta'}^{-1}\in\mathfrak{A}_Q$, for all $\Theta'\in\mathfrak{I}_{P,Q}$ and, therefore, $|\mathfrak{I}_{P,Q}|\leq |\mathfrak{A}_Q|$. The opposite inequality is also verified, because, given $\Theta\in\mathfrak{I}_{P,Q}$, it is $\Theta'\Theta\in\mathfrak{I}_{P,Q}$, for all $\Theta'\in\mathfrak{A}_Q$. The equality also implies that the last assertion of the result is then an immediate consequence of the first part of the reasoning. \hfill $\Box$

\vspace{0.5cm}

Hereafter, we focus our study on the values $\Delta_s(z)$. The values $\Delta_{[P]}(z)$ needs a comprehensive analysis of the isotopic classes of partial Latin squares and will be considered in a further study. Firstly, it raises the natural question of whether it is possible to obtain some general expression which determines these values for some specific size or cycle structure. So, for instance, it is immediate to see that $\Delta_s((1^n,1^n,1^n))=|\mathcal{PLS}_{n,s}|$ and, since $\mathcal{PLS}_{n,n^2}=\mathcal{LS}_n$, it is also clear that $\Delta_{n^2}(z)=\Delta(z)$. In this regard, let us study some cases in which a general formula is given:

\begin{prp}\label{prpD2} Let $s\in [n^2]$. It is verified that:
$$\Delta_s((n,n,1^n))=\begin{cases}\frac {{n!}^2}{k! \cdot {(n-k)!}^2}, \text { if } \exists k\in [n] \text{ s.t. } s=k\cdot n,\\0, \text{ otherwise}.\end{cases}$$
\end{prp}

{\bf Proof.}  Let $\Theta=(\alpha,\beta,\mathrm{Id})\in\mathfrak{I}_{(n,n,1^n)}$, where ${\mathrm{Id}}$ denotes the trivial permutation, and $P\in\mathcal{PLS}_{\Theta,s}$. Since the $\Theta$-decomposition of $P$ is only formed by $P$ itself and $LMC_z=\{(n,n)\}$, then Lemma \ref{lemR1} implies $s=k\cdot n$, for some $k\in [n]$. Thus, $O(P)$ is decomposed under the action of $\Theta$ into $k$ orbits of length $n$. Specifically, there exist exactly $k$ distinct columns $c_1,c_2,...,c_k\in [n]$ and $k$ distinct symbols $s_1,s_2,...,s_k\in [n]$, such that $(1,c_i,s_i)\in O(P)$, for all $i\in [k]$. The $k$ orbits of $O(P)$ under $\Theta$ are then the sets $\{(\alpha^t(1),\beta^t(c_i),s_i)\,\mid\, t\in [n]\}$, with $i\in [k]$.

Every element of $\mathcal{PLS}_{\Theta,k\cdot n}$ is therefore univocally determined by the choice of the columns $c_i$ and symbols $s_i$. Namely, there exist $\left(\begin{array}{c}n\\k\end{array}\right)$ possible ways of choosing the $k$ columns and, once they have been selected, there exist $\frac {n!}{(n-k)!}$ different ways of assigning $k$ symbols to the cells $(1,c_1),...,(1,c_n)$. So, $\Delta_s(z)=\left(\begin{array}{c}n\\k\end{array}\right)\cdot \frac {n!}{(n-k)!}=\frac {{n!}^2}{k! \cdot {(n-k)!}^2}$.\hfill $\Box$

\vspace{0.5cm}

\begin{prp}\label{prpD3} It is verified that $\Delta_n((n,n,n))=n^2$. Besides, if $n>2$:
$$\Delta_{2n}((n,n,n))=\frac {n^2\cdot (n-1)\cdot (n-2)}2.$$
\end{prp}

{\bf Proof.}  Let $\Theta=(\alpha,\beta,\gamma)\in\mathfrak{I}_{(n,n,n)}$. Similarly to the proof of Proposition \ref{prpD2}, every partial Latin square of $\mathcal{PLS}_{\Theta,n}$ is univocally determined by the only non-empty cell $(1,c)$ of its first row. Specifically, $c$ can be selected from $n$ different columns and the symbol of $(1,c)$ will be able to be chosen between $n$ candidates. So, $\Delta_n((n,n,n))=n^2$.

Analogously, if $n>2$, every partial Latin square of $\mathcal{PLS}_{\Theta,2n}$ is univocally determined by the only two non-empty cells of its first row. Their corresponding columns $c_1, c_2\in [n]$ can be selected of $\left(\begin{array}{c}n\\2\end{array}\right)$ different ways. Now, the symbol $s_1$ of the cell $(1,c_1)$ can be chosen from $n$ possibilities, but, once it has been chosen, there exist only $n-2$ candidates for the symbol $s_2$ of the cell $(1,c_2)$. To see it, let $t\in [n-1]$ be such that $\beta^t(c_1)=c_2$. So, the triple $(\alpha^t(1),c_2,\gamma^t(s_1))$ belongs to $O(P)$. Moreover, since $z_{\gamma}=n$, then $s_1\neq \gamma^t(s_1)$ and the Latin square condition implies $s_2\not\in\{s_1,\gamma^t(s_1)\}$. Thus, $\Delta_{2n}((n,n,n))=\left(\begin{array}{c}n\\2\end{array}\right)\cdot n\cdot (n-2)=\frac {n^2\cdot (n-1)\cdot (n-2)}2$. \hfill $\Box$

\vspace{0.5cm}

\begin{thm}\label{thmD1} Given $z=(z_1,z_2,z_3)\in {\mathcal{CS}_{\mathfrak{A}_{\mathcal{P}}}}_n$:
$$\Delta_{\mathfrak{l}_z}(z)=\sum_{\scriptsize \begin{array}{c}(i,j)\in \mathrm{LCM}_z\\\text{s.t. } \mathrm{lcm}(i,j)=\mathfrak{l}_z\end{array}} z_{1i}\cdot z_{2j}\cdot \mathrm{gcd}(i,j)\cdot \sum_{\scriptsize \begin{array}{c}k\in [n]\\\text{s.t. } (i,j,k)\in \mathrm{LCM}_n\end{array}}k\cdot z_{3k}.$$
\end{thm}

{\bf Proof.} Given $\Theta=(\alpha,\beta,\gamma)\in\mathfrak{I}_z$, let $P\in \mathcal{PLS}_{\Theta}$ be such that $|P|=\mathfrak{l}_z$. From Lemma \ref{lemR1} and Proposition \ref{prpR1}, there must exist only one non-empty block $B$ in the $\Theta$-decomposition of $P$. Specifically, $B$ must be an $i\times j$-block of size $\mathrm{lcm}(i,j)=\mathfrak{l}_z$, where $(i,j)\in \mathrm{LCM}_z$. There exist $z_{1i}\cdot z_{2j}$ possible blocks in this way.

Moreover, $O(B)$ must be composed by all the triples of one of the $\mathrm{gcd}(i,j)$ orbits induced on $B$ by the action of $\Theta$. If $(r,c,s)\in [n]^3$ is one of these triples, then the symbol $s$ must be one of the $k\cdot z_{3k}$ elements of a $k$-cycle of $\gamma$ such that $(i,j,k)\in \mathrm{LCM}_n$. The result follows then by considering all the previous possibilities. \hfill $\Box$

\vspace{0.5cm}

\begin{cor}\label{corD1} Let $P\in \mathcal{PLS}_{n,1}$. Given $z=(z_1,z_2,$ $z_3)\in {\mathcal{CS}_{\mathfrak{A}_{\mathcal{P}}}}_n$:
$$\Delta_{[P]}(z)=\Delta_1(z)=z_{11}\cdot z_{21}\cdot z_{31}.$$
\end{cor}

{\bf Proof.} Since there exists only one isotopic class of partial Latin squares of size 1, it is $\Delta_{[P]}(z)=\Delta_1(z)$. Now, if $\mathfrak{l}_z>1$, then $(1,1)\not\in \mathrm{LCM}_z$. So, $z_{11}\cdot z_{21}\cdot z_{31}=0$ and the result holds. Finally, if $\mathfrak{l}_z=1$, then it is enough to observe that it must be $(i,j,k)=(1,1,1)$ in the formula of Theorem \ref{thmD1}. \hfill $\Box$

\vspace{0.5cm}

The number $\Delta_s(z)$ can also obviously be obtained if the set $\mathcal{PLS}_{\Theta,n}$ is known for some $\Theta\in\mathfrak{I}_z$. In order to determine this set, let us observe that, analogously to the case of Latin squares \cite{Euler86}, $\mathcal{PLS}_n$ can be identified \cite{Kumar99} with the set $\mathcal{FS}_{\mathcal{P}_n}$ of feasible solutions of the integer program:
$$\begin{array}{l}\min \sum_{r, c, s\in [n]} w_{rcs}\cdot x_{rcs},\\
\begin{array}{ll}\text{subject to } & \sum_{r\in [n]}x_{rcs}\leq 1, \forall c, s\in [n],\\
         \ & \sum_{c\in [n]}x_{rcs}\leq 1, \forall r, s\in [n],\\
         \ & \sum_{s\in [n]}x_{rcs}\leq 1, \forall r, c\in [n],\\
         \ & x_{rcs}\in \{0,1\}, \forall r, c, s \in [n],\end{array}\end{array} \hspace{2cm} (1)$$
where $w_{rcs}$ are real weights for all $r, c, s \in [n]$. Specifically, it is enough to define the map $\varphi_n:\mathcal{PLS}_n\rightarrow \mathcal{FS}_{\mathcal{P}_n}$, such that, given $P\in\mathcal{PLS}_n$, it is $\varphi_n(P)=(x^P_{111},...,x^P_{11n},x^P_{121},...,x^P_{nnn})$, where, $x^P_{rcs}=1$ if $(r,c,s)\in O(P)$ and $0$, otherwise. The restriction of $\varphi_n$ to $\mathcal{PLS}_{\Theta}$ and $\mathcal{PLS}_{\Theta,m}$ assures the truthfulness of the following result:

\begin{prp} \label{prpP1} Given $\Theta=(\alpha,\beta,\gamma)\in \mathfrak{I}_n$, there exists a bijection between $\mathcal{PLS}_{\Theta}$ and the set of feasible solutions of the equation system which results after adding to $(1)$ the constraints:
$$x_{rcs}=x_{\alpha(r)\beta(c)\gamma(s)},\forall r, c, s \in [n].$$
Moreover, given $m\in [n^2]$, if the equation:
$$\sum_{r,c,s\in [n]} x_{rcs} = m$$
is also added, then there exists a bijection between $\mathcal{PLS}_{\Theta,m}$ and the set of feasible solutions of the resulting equation system. \hfill $\Box$
\end{prp}

\vspace{0.5cm}

Proposition \ref{prpP1} implies $\mathcal{PLS}_{\Theta,m}$ to be determined by $2n^3+3n^2+1$ polynomial equations of degree $1$ and $2$ in $n^3$ variables:

\begin{cor} \label{corP1} Given $\Theta=(\alpha,\beta,\gamma)\in \mathfrak{I}_n$ and $m\in[n^2]$, $\mathcal{PLS}_{\Theta,m}$ is the set of zeros of the ideal $I=\langle\,(\sum_{r\in [n]}x_{rcs})\cdot (1-\sum_{r\in [n]}x_{rcs})=0\,\mid\, c, s\in [n]\,\rangle + \langle\,(\sum_{c\in [n]}x_{rcs})\cdot (1-\sum_{c\in [n]}x_{rcs})\,\mid\, r, s\in [n]\,\rangle + \langle\,(\sum_{s\in [n]}x_{rcs})\cdot (1-\sum_{s\in [n]}x_{rcs})\,\mid\, r, c\in [n]\,\rangle + \langle\,x_{rcs}\cdot\left(1-x_{rcs}\right)\,\mid\, r, c, s \in [n]\,\rangle + \langle\,x_{rcs}-x_{\alpha(r)\beta(c)\gamma(s)}\,\mid\, r, c, s \in [n] \,\rangle + \langle\,m-\sum_{r,c,s\in [n]} x_{rcs}\,\rangle\,\subseteq \mathbb{Q}[{\bf x_n}] = \mathbb{Q}[x_{111},...,$ $x_{nnn}]$. \hfill $\Box$
\end{cor}

\vspace{0.5cm}

The ideal $I$ of Corollary \ref{corP1} is {\em zero-dimensional}, i.e., there exists only a finite number of solutions of the corresponding system of polynomial equations. Moreover, $I\cap \mathbb{Q}[x_{rcs}]=\langle\,x_{rcs}\cdot(1-x_{rcs})\,\rangle\subseteq I$, for all $r,c,s\in [n]$, and, therefore, Proposition 2.7 of \cite{Cox98} implies $I$ to be {\em radical}, i.e., any polynomial $p({\bf x_n})$ belongs to $I$ whenever there exists $t\in \mathbb{N}$ such that $p({\bf x_n})^t\in I$. Since the affine variety defined by $I$ is $V(I)=\mathcal{PLS}_{\Theta,m}$, then Theorem 2.10 of \cite{Cox98} assures $\Delta_m(z_{\Theta})=|V(I)|=\mathrm{dim}_{\mathbb{Q}}(\mathbb{Q}[{\bf x_n}]/I)$, which can be computed from any Gr\"obner basis of $I$, with respect to any term ordering. So, for instance, {\sc Singular} \cite{Decker11} has been used in order to obtain the values of $\Delta_s(z)$ and $\Delta_{\mathcal{P}}(z)$ for each parastrophic class of ${\mathcal{CS}_{\mathfrak{A}_{\mathcal{P}}}}_n$, where $n\leq 4$. These values are shown in Tables 2 and 3, where the blank cells correspond to those $s\not\in\mathrm{Sizes}(z)$.

\begin{table}
\begin{center}{\footnotesize
\begin{tabular}{|c|c|c|c|c|c|c|c|c|c|c|c|} \hline
\multirow{3}{*}{$n$} & \multirow{3}{*}{$z$} & \multicolumn{9}{|c|}{$\Delta_s(z)$}& \multirow{3}{*}{$\Delta_{\mathcal{P}}(z)$}\\ \cline{3-11}
\ & \ & \multicolumn{9}{|c|}{$s$}&\ \\ \cline{3-11}
\ & \  & 1 & 2 & 3 & 4&5&6&7&8&9&\ \\ \hline
1 & $(1,1,1)$ & 1 & &  &  &  & &  & &  &1\\ \hline
\multirow{3}{*}{2} & (2,2,2) &   & 4 &   & 0& &  &  &  & &4\\ \cline{2-12}
\ & $(2,2,1^2)$ &   & 4 &   & 2& &  &  &  & &6\\ \cline{2-12}
\ & $(1^2,1^2,1^2)$ & 8 & 16 & 8 & 2& &  &  &  & &34\\ \hline
\multirow{7}{*}{3} & (3,3,3) &  &  &9 & &  &9 &  &  &3&21\\\cline{2-12}
\ & (3,3,2$\cdot$1) &  &  &3 &  &  &  &  &  & &3\\\cline{2-12}
\ &$(3,3,1^3)$ &  &  &9 &  &  &18 &  &  &6&33\\\cline{2-12}
\ &(2$\cdot$1,2$\cdot$1,2$\cdot$1) &1 &10 &10 &24 &24 &20 &20 &4 &4&117\\\cline{2-12}
\ &(2$\cdot$1,2$\cdot$1,$1^3$) &3 &6 &18 &6 &18 &  &  &  & &51\\\cline{2-12}
\ &(2$\cdot$1,$1^3$,$1^3$) &9 &18 &6 &  &  &  &  &  & &33\\\cline{2-12}
\ &$(1^3,1^3,1^3)$ &27 &270 &1278 &3078 &3834 &2412 &756 &108 &12&11775 \\\hline
\end{tabular}}
\end{center}
\caption{$\Delta_s(z)$ and $\Delta_{\mathcal{P}}(z)$ for each parastrophic class of ${\mathcal{CS}_{\mathfrak{A}_{\mathcal{P}}}}_n$, where $n\leq 3$.}
\end{table}

\renewcommand
{\tabcolsep}{1pt}

\begin{sidewaystable}
\begin{center}{\scriptsize
\begin{tabular}{|c|c|c|c|c|c|c|c|c|c|c|c|c|c|c|c|c|c|c|} \hline
\multirow{3}{*}{$z$} & \multicolumn{16}{|c|}{$\Delta_s(z)$}& \multirow{3}{*}{$\Delta_{\mathcal{P}}(z)$}\\ \cline{2-17}
\ & \multicolumn{16}{|c|}{$s$}&\ \\ \cline{2-17}
\  & 1 & 2 & 3 & 4&5&6&7&8&9&10&11&12&13&14&15&16& \\ \hline
(4,4,4) &  & &  &16  &  & &  &48  & &  &  &32  & &  & &0 &96 \\ \hline
(4,4,3$\cdot$1) &  & &  &4  &  & &  &  & &  &  &  & &  & & &4 \\ \hline
(4,4,$2^2$) &  & &  &16  &  & &  &56  & &  &  &32  & &  & &8 &112 \\ \hline
(4,4,2$\cdot 1^2$) &  & &  &16  &  & &  &64  & &  &  &64  & &  & &8 &152 \\ \hline
(4,4,$1^4$) &  & &  &16  &  & &  &72  & &  &  &96  & &  & &24 &208 \\\hline
(3$\cdot$1,3$\cdot$1,3$\cdot$1) &1  & &18  &18  &  &90 &90  &  &165 &165  &  &99  &99 &  &9 &9 &763 \\\hline
(3$\cdot$1,3$\cdot$1,2$\cdot1^2$) &2  & &6  &12  &  &6 &12  &  & &  &  &  & &  & & &38 \\\hline
(3$\cdot$1,3$\cdot$1,$1^4$) &4  & &12  &48  &  &36 &144  &  &24 &96  &  &  & &  & & &364 \\\hline
(3$\cdot$1,$2^2$,$2^2$) &  &8 &  &8  &  & &  &  & &  &  &  & &  & & &16 \\\hline
(3$\cdot$1,$2^2$,2$\cdot1^2$) &  &4 &  &  &  & &  &  & &  &  &  & &  & & &4 \\\hline
(3$\cdot$1,2$\cdot1^2$,2$\cdot1^2$) &4  &4 &8  &4  &  & &  &  & &  &  &  & &  & & &20  \\\hline
(3$\cdot$1,2$\cdot1^2$,$1^4$) &8  &12 &  &  &  & &  &  & &  &  &  & &  & & &20  \\\hline
(3$\cdot$1,$1^4$,$1^4$) &16  &72 &96  &24  &  & &  &  & &  &  &  & &  & & &208  \\\hline
($2^2$,$2^2$,$2^2$) &  &32 &  &352  &  &1664 &  &3552  & &3328  &  &1408  & &256  & &32 &10624  \\\hline
($2^2$,$2^2$,2$\cdot1^2$)&  &32 &  &360 &  &1792  &  &4152  & &4416  &  &2048  & &384  & &32 &13216  \\\hline
($2^2$,$2^2$,$1^4$)&  &32 &  &368  &  &1920 &  &4800  & &5760  &  &3264  & &768  & &96 &17008  \\\hline
($2^2$,2$\cdot1^2$,2$\cdot1^2$)&  &24 &  &192  &  &640 &  &880  & &416  &  &32  & &  & & &2184  \\\hline
($2^2$,2$\cdot1^2$,$1^4$)&  &16 &  &72  &  &96 &  &24  & &  &  &  & &  & & &208  \\\hline
(2$\cdot1^2$,2$\cdot1^2$,2$\cdot1^2$)&8  &32 &136  &336  &752  &1440 &1904  &2856  &2400 &2608  &1504  &1056  &448 &224  &64 &16 &15784\\\hline
(2$\cdot1^2$,2$\cdot1^2$,$1^4$)&16  &88 &272  &736  &1344  &1632 &1728  &1008  & &  &  &  & &  & & &6824  \\\hline
(2$\cdot1^2$,$1^4$,$1^4$)&28  &352 &2208  &6504  &9792  &7104 &2112  &216  & &  &  &  & &  & & &28352\\\hline
($1^4$,$1^4$,$1^4$)&64  &1728 &25920  &239760  &1437696  &5728896 &15326208  &27534816  &32971008 &25941504  &13153536  &4215744  &847872 &110592  &9216 &576 &127545136\\\hline
\end{tabular}}
\end{center}
\caption{$\Delta_s(z)$ and $\Delta_{\mathcal{P}}(z)$ for each parastrophic class of ${\mathcal{CS}_{\mathfrak{A}_{\mathcal{P}}}}_4$.}
\end{sidewaystable}

\section{$\Theta$-completable partial Latin squares.}

In the Introduction, given $z\in {\mathcal{CS}_{\mathfrak{A}_{\mathcal{P}}}}_n$ and $\Theta\in\mathfrak{I}_z$, it has been indicated that a partial Latin square $P\in\mathcal{PLS}_{\Theta}$ can be used in the computation of $\Delta(z)$, if a multiplicative factor ($P$-coefficient of symmetry of $\Theta$ \cite{FalconMartinJSC07}) $c_P\in\mathbb{N}$ is found such that $\Delta(z)=c_P\cdot |\mathcal{LS}_{\Theta,P}|$. In this regard, let us finish the present study with a theoretical basis for this concept of coefficient of symmetry. To do it, it is necessary to generalize the traditional concept of completability. Specifically, $P$ will be said to be {\em $\Theta$-completable} if $\mathcal{LS}_{\Theta,P}\neq\emptyset$. So, the traditional completability corresponds to the trivial isotopism $\Theta=(\mathrm{Id},\mathrm{Id},\mathrm{Id})$. Moreover, let us observe that, if a partial Latin square is $\Theta$-completable, then it is also completable in the traditional way.

It can be easily checked that every partial Latin square of order $n\leq 2$ is $\Theta$-completable whenever $\Theta$ is non-trivial. An example of non-trivial isotopism for which there exists a related partial Latin square which is neither $\Theta$-completable nor traditionally completable is $\Theta=((12)(3),(12)(3),(12)(3))\in\mathfrak{A}_3$. A partial Latin square in such conditions is:{\scriptsize
$$\left(\begin{array}{ccc}3 & \cdot & 2\\
                          \cdot & 3 & 1\\
                          2 & 1 & \cdot\end{array}\right).$$}

An example where it is possible to observe the difference between both concepts is given if $\Theta=((12)(34),(12)(34),(12)(3)(4))\in\mathfrak{A}_4$. In this case, the following partial Latin square is not $\Theta$-completable, but it is completable in the traditional way:{\scriptsize
$$\left(\begin{array}{cccc}3 & 4 & \cdot & \cdot\\
                           4 & 3 & \cdot & \cdot\\
                           \cdot & \cdot & \cdot & \cdot\\
                           \cdot & \cdot & \cdot & \cdot\end{array}\right).$$}

Let $C_{\Theta}$ denote the set of $\Theta$-completable partial Latin squares and let $C_{\Theta,s}=C_{\Theta}\cap\mathcal{PLS}_{\Theta,s}$. The cardinalities of these sets only depends on the parastrophic class of the cycle structure of $\Theta$:

\begin{lem}\label{lemT1} Let $\Theta_1,\Theta_2\in \mathfrak{I}_n$ be such that $[z_{\Theta_1}]=[z_{\Theta_2}]$. Then, $|C_{\Theta_1,s}|=|C_{\Theta_2,s}|$, for all $s\in[n^2]$. As a consequence, $|C_{\Theta_1}|=|C_{\Theta_2}|$.
\end{lem}

{\bf Proof.} Let $\pi\in S_3$ be such that $z_{\Theta_1}=z_{\Theta^{\pi}_2}$. Given $s\in [n^2]$ and $P\in C_{\Theta_1,s}$, there exists $L\in \mathcal{LS}_{\Theta_1}$ such that $O(P)\subseteq O(L)$. Besides, since $\Theta_1$ and $\Theta^{\pi}_2$ are conjugate, then there exists $\Theta\in\mathfrak{I}_n$ such that $\Theta^{\pi}_2=\Theta\Theta_1\Theta^{-1}$. Thus, $(P^{\Theta})^{\pi^{-1}}\in\mathcal{PLS}_{\Theta_2}$, $(L^{\Theta})^{\pi^{-1}}\in\mathcal{LS}_{\Theta_2}$ and $O((P^{\Theta})^{\pi^{-1}})\subseteq O((L^{\Theta})^{\pi^{-1}})$. Since $|(P^{\Theta})^{\pi^{-1}}|=|P|$, then $|C_{\Theta_1,s}|\leq |C_{\Theta_2,s}|$. The opposite inequality is analogously proven and the consequence is immediate, because $|C_{\Theta_1}|=\sum_{s\in [n^2]}|C_{\Theta_1,s}|=\sum_{s\in [n^2]}|C_{\Theta_2,s}|=|C_{\Theta_2}|$. \hfill $\Box$

\vspace{0.5cm}

From the previous result, it is natural to define the numbers $\mathfrak{c}_z$ and $\mathfrak{c}_{z,s}$ as the respective cardinalities of $C_{\Theta}$ and $C_{\Theta,s}$, for any $\Theta\in\mathfrak{I}_z$. The following result holds:

\begin{thm}\label{thmT1} Let $z\in {\mathcal{CS}_{\mathfrak{A}_{\mathcal{P}}}}_n$ and $\Theta\in\mathfrak{I}_z$. It is verified that $P\in C_{\Theta}$ if and only if $\mathcal{PLS}_{\Theta,[P]}\subseteq C_{\Theta,|P|}$. As a consequence:
$$\mathfrak{c}_{z,s}=\sum_{\scriptsize \begin{array}{c}[P]\in \mathcal{PLS}_{\Theta,s}/\sim\\ \text{s.t. }[P]\cap C_{\Theta}\neq \emptyset\end{array}}
\Delta_{[P]}(z).$$
\end{thm}

{\bf Proof.} The sufficient condition is immediate. So, let us consider $P\in C_{\Theta}$ and $Q\in \mathcal{PLS}_{\Theta,[P]}$. Since $Q\in [P]$ and $\Theta\in\mathfrak{A}_P\cap\mathfrak{A}_Q$, Theorem \ref{thmD2} implies that there exists $\Theta'\in\mathfrak{I}_{P,Q}$ such that $\Theta=\Theta'\Theta\Theta'^{-1}$. Thus, $\Theta\Theta'=\Theta'\Theta$. Now, let $L\in\mathcal{LS}_{\Theta}$ be such that $O(P)\subseteq O(L)$. It must be then  $L^{\Theta'}\in\mathcal{LS}_{\Theta}$, because $(L^{\Theta'})^{\Theta}=(L^{\Theta})^{\Theta'}=L^{\Theta'}$. Since $|Q|=|P|$ and $O(Q)=O(P^{\Theta'})\subseteq O(L^{\Theta'})$, then $Q\in C_{\Theta,|P|}$ and the first claim is verified. The consequence is then immediate. \hfill $\Box$

\vspace{0.5cm}

The previous theorem implies that it is enough to check the completability of one element of each isotopic class of partial Latin squares.
Moreover, it is convenient to do it in increasing order of the size, because, given $\Theta\in\mathfrak{I}_n$ and $P,Q\in\mathcal{PLS}_{\Theta}$ such that $O(P)\subseteq O(Q)$, if $P$ is not $\Theta$-completable, neither is $Q$. Taking into account this strategy, the numbers $\mathfrak{c}_{z,s}$ and $\mathfrak{c}_z$ have been obtained (Table 5) for each non-trivial parastrophic class of $\mathcal{CS}_{{\mathfrak{A}}_n}$, where $n\leq 4$.

\renewcommand
{\tabcolsep}{2pt}

\begin{table}
\begin{center}{\scriptsize
\begin{tabular}{|c|c|c|c|c|c|c|c|c|c|c|c|c|c|c|c|c|c|c|} \hline
\multirow{3}{*}{$n$} & \multirow{3}{*}{$z$} & \multicolumn{16}{|c|}{$\mathfrak{c}_{z,s}$}& \multirow{3}{*}{$\mathfrak{c}_z$}\\ \cline{3-18}
\ & \ & \multicolumn{16}{|c|}{$s$}&\ \\ \cline{3-18}
\ & \  & 1 & 2 & 3 & 4&5&6&7&8&9&10&11&12&13&14&15&16&\ \\ \hline
1 & $(1,1,1)$ & 1 & &  &  &  & &  & &  & & & & & & & &1\\ \hline
2 & $(2,2,1^2)$ &   & 4 &   & 2& &  &  &  & & & & & & & & &6\\ \hline
\multirow{3}{*}{3} & (3,3,3) &  &  &9 & &  &9 &  &  &3& & & & & & & &21\\\cline{2-19}
\ &$(3,3,1^3)$ &  &  &9 &  &  &18 &  &  &6& & & & & & & &33\\\cline{2-19}
\ &(2$\cdot$1,2$\cdot$1,2$\cdot$1) &1 &10 &10 &24 &24 &16 &16 &4 &4& & & & & & & &109\\\hline
\multirow{8}{*}{$4$} & (4,4,$2^2$) &  & &  &16  &  & &  &40  & &  &  &32  & &  & &8 &96 \\ \cline{2-19}
\ &(4,4,2$\cdot 1^2$) &  & &  &16  &  & &  &40  & &  &  &32  & &  & &8 &96 \\ \cline{2-19}
\ &(4,4,$1^4$) &  & &  &16  &  & &  &72  & &  &  &96  & &  & &24 &208 \\\cline{2-19}
\ &(3$\cdot$1,3$\cdot$1,3$\cdot$1) &1  & &18  &18  &  &90 &90  &  &90 &90  &  &45  &45 &  &9 &9 &505 \\\cline{2-19}
\ &($2^2$,$2^2$,$2^2$) &  &32 &  &352  &  &1408 &  &2144  & &1792  &  &896  & &256  & &32 &6912  \\\cline{2-19}
\ &($2^2$,$2^2$,2$\cdot1^2$)&  &32 &  &336 &  &1344  &  &2144  & &1792  &  &896  & &256  & &32 &6832  \\\cline{2-19}
\ &($2^2$,$2^2$,$1^4$)&  &32 &  &368  &  &1728 &  &3792  & &4224  &  &2496  & &768  & &96 &13504  \\\cline{2-19}
\ &(2$\cdot1^2$,2$\cdot1^2$,2$\cdot1^2$)&8  &32 &136  &200  &528  &784 &1328  &1560  &1760 &1568  &1248  &800  &448 &192  &64 &16 &10672\\\hline
\end{tabular}}
\end{center}
\caption{$\mathfrak{c}_{z,s}$ and $\mathfrak{c}_{z}$ for each non-trivial parastrophic class of ${\mathcal{CS}_{\mathfrak{A}_{\mathcal{P}}}}_n$, where $n\leq 4$.}
\end{table}

\vspace{0.5cm}

Given $z\in \mathcal{CS}_{\mathfrak{A}_{\mathcal{P}_n}}$ and $\Theta\in\mathfrak{I}_z$, a set $\{P_1,P_2,...,P_m\}$ of $\Theta$-completable partial Latin squares will be said to be a {\em basis} of $\mathcal{LS}_{\Theta}$ if $\bigcup_{i\in [m]} \mathcal{LS}_{\Theta,P_i} = \mathcal{LS}_{\Theta}$ and $\mathcal{LS}_{\Theta,P_i}\cap\mathcal{LS}_{\Theta,P_j}=\emptyset$, whenever $i\neq j$. In this case, $\Delta(z)=\sum_{i\in [m]}|\mathcal{LS}_{\Theta,P_i}|$. Let us observe that, from a computational point of view, it is interesting to determine a basis of $\mathcal{LS}_{\Theta}$ such that the sizes of its elements are as great as it is possible, because then, for each $P_i$, it would be feasible to add to the constraints of Proposition \ref{prpP1}, all those of the form $x_{rcs}=1$, if $(r,c,s)\in O(P_i)$. The calculus of the corresponding Gr\"obner basis would be then more efficient and it would allow to obtain new values $\Delta(z)$. The following result is proven:

\begin{lem}\label{lemT2} Let $S\subseteq [n]^2$ and $\Theta=(\alpha,\beta,\gamma)\in\mathfrak{A}_n$. Each of the following sets is non-empty if and only if it is a basis of $\mathcal{LS}_{\Theta}$:
$$S_{RC}=\{P\in C_{\Theta}\,\mid\, (r,c,s)\in O(P) \Leftrightarrow (r,c)\in S\},$$
$$S_{RS}=\{P\in C_{\Theta}\,\mid\, (r,c,s)\in O(P) \Leftrightarrow (r,s)\in S\},$$
$$S_{CS}=\{P\in C_{\Theta}\,\mid\, (r,c,s)\in O(P) \Leftrightarrow (c,s)\in S\}.$$ \hfill $\Box$
\end{lem}

{\bf Proof.} The sufficient condition is immediate. In order to see the necessary condition, let us prove that $S_{RC}$ is a basis of $\mathcal{LS}_{\Theta}$; the other cases are similar. Since $S_{RC}$ is non-empty, there exists $P_0\in C_{\Theta}$ such that $(r,c,s)\in O(P_0) \Leftrightarrow (r,c)\in S$. Given $L\in \mathcal{LS}_{\Theta}$, let $P\in\mathcal{PLS}_n$ be such that $O(P)=\{(r,c,s)\in O(L)\,\mid\, (r,c)\in S\}$ and let us consider $(r,c,s)\in O(P)$. It must be $(r,c)\in S$ and so, there must exist $s_0\in [n]$ such that $(r,c,s_0)\in O(P_0)$. Since $P_0\in \mathcal{PLS}_{\Theta}$, it must be $(\alpha(r),\beta(c),\gamma(s_0))\in O(P_0)$ and therefore, $(\alpha(r),\beta(c))\in S$. Hence, $(\alpha(r),\beta(c),\gamma(s))\in O(P)$ and thus, $P\in\mathcal{PLS}_{\Theta}$. As a consequence, $P$ is $\Theta$-completable and then,
$\mathcal{LS}_{\Theta}\subseteq \bigcup_{Q\in S_{RC}} \mathcal{LS}_{\Theta,Q}$. Indeed, both sets are equal because $\mathcal{LS}_{\Theta,Q}\subseteq\mathcal{LS}_{\Theta}$, for all $Q\in S_{RC}$. Finally, given two distinct elements $Q,Q'\in S_{RC}$, it must exist $(r,c)\in S$ and $s\in [n]$ such that $(r,c,s)\in O(Q)\setminus O(Q')$. It implies that $\mathcal{LS}_{\Theta,Q}\cap\mathcal{LS}_{\Theta,Q'}=\emptyset$ and therefore, $S_{RC}$ is a basis of $\mathcal{LS}_{\Theta}$. \hfill $\Box$

\vspace{0.5cm}

A special case appears when $|\mathcal{LS}_{\Theta,P_i}|=|\mathcal{LS}_{\Theta,P_j}|$, for all $i,j\in [m]$. Such a basis will be called {\em homogeneous} and it will be verified that $\Delta(z)=m\cdot |\mathcal{LS}_{\Theta,P_i}|$, for all $i\in[m]$. The cardinality $m$ of the homogeneous basis would be therefore the $P_i$-coefficient of symmetry of $\Theta$, for all $i\in [m]$. Although a comprehensive study must be developed in this regard, let us finish the current paper with a result with gives a theoretical support to the majority of the coefficients of symmetry which were used in \cite{FalconMartinJSC07}:

\begin{thm}\label{thmT2} Let $z=(z_1,z_2,z_3)\in \mathcal{CS}_{\mathfrak{A}_n}$ be such that $z_{11}\cdot z_{21}\cdot z_{31}\neq 0$. Let $\Theta=(\alpha,\beta,\gamma)\in\mathfrak{I}_z$ and $S=\{(i,j)\in [n]^2\,\mid\, i\in\alpha_{\infty}, j\in \beta_{\infty}\}$. It is verified that $S_{RC}$ is an homogeneous basis of $\mathcal{LS}_{\Theta}$ of cardinality $|\mathcal{LS}_{z_{11}}|$.
\end{thm}

{\bf Proof.} From the hypothesis, it must be $z_1=z_2=z_3$ (\cite{McKay05}, Theorem 1). Furthermore, given $P\in\mathcal{PLS}_{\Theta}$, the corresponding block $P_{\infty\infty}$ of the $\Theta$-decomposition of $P$ is a $z_{11}\times z_{11}$-array, such that each of its non-filled cells must contain one of the $z_{11}$ fixed symbols of $\gamma$, i.e., it is a Latin subsquare of $P$ of order $z_{11}$. Thus, since $\Theta\in \mathfrak{A}_n$, Lemma \ref{lemT2} implies the set $S_{RC}$ to be a basis of $\mathcal{LS}_{\Theta}$ of $|\mathcal{LS}_{z_{11}}|$ elements. Now, let us consider two distinct elements $Q,Q'\in S_{RC}$. Given $L\in\mathcal{LS}_{\Theta,Q}$, let us define the Latin square $L'\in\mathcal{LS}_n$ such that $O(L')=\{(r,c,s)\in [n]^3\,\mid\, (r,c,s)\in O(Q') \text{ if } (r,c)\in S, \text{ or } (r,c,s)\in O(L), \text{ otherwise}\}$, i.e., the only difference of $L'$ with respect to $L$ is the block $L'_{\infty\infty}$, which is $Q'$ instead of $Q$. Since $L\in \mathcal{LS}_{\Theta}$ and $Q'\in\mathcal{PLS}_{\Theta}$, it must be $L'\in\mathcal{LS}_{\Theta}$. Hence, $|\mathcal{LS}_{\Theta,Q}|\leq |\mathcal{LS}_{\Theta,Q'}|$. The opposite inequality is analogously proven and, therefore, $S_{RC}$ is homogeneous. \hfill $\Box$

\section{Final remarks and further work.}

In the current paper, it has been dealt with the set of autotopisms of partial Latin squares in order to develop further techniques which allow to improve some results about the set of autotopisms of Latin squares, such as those related with the obtention of the values $\Delta(z)$. In Section 2, the cardinality of ${\mathcal{CS}_{\mathfrak{A}_{\mathcal{P}}}}_n$ has been studied and a lower bound has been determined. Although it can be obtained by an exhaustive search once Lemma \ref{lem0} is implemented in a computer procedure, it raises the question of whether it is possible to obtain a general formula for $|{\mathcal{CS}_{\mathfrak{A}_{\mathcal{P}}}}_n|$. A similar question appears in Section 4 with the values $\Delta_{[P]}(z)$, for which a comprehensive study of isotopic classes of $\mathcal{PLS}_n$ would be necessary. It would also be useful in order to improve the computation and increase the order $n\leq 4$ which have been used in the examples of the present paper. Finally, once a theoretical basis has been exposed in Section 5 for the concept of coefficient of symmetry of an autotopism, it seems that an exhaustive study in this regard would be necessary to solve some of the problems of computation related to the calculus of the values $\Delta(z)$.


\begin{thebibliography}{00}

\bibitem{Bryant09} D. Bryant, M. Buchanan and I. M. Wanless, The spectrum for quasigroups with cyclic automorphisms and additional symmetries, Discrete Math. 304 (2009), 821-833.

\bibitem{Cox98} D. A. Cox, J. B. Little and D. O'Shea, Using Algebraic Geometry, Springer-Verlag, New York, 1998.

\bibitem{Decker11} W. Decker, G.-M. Greuel, G. Pfister and H. Sch\"onemann, {\sc Singular} 3-1-3. A computer algebra system for polynomial computations, 2011. http://www.singular.uni-kl.de.

\bibitem{Euler86} R. Euler, R. E. Burkard and R. Grommes, On Latin squares and the facial structure of related polytopes, Discrete Math. 62 (1986), 155-181.

\bibitem{Falcon06} R. M. Falc\'on, Latin squares associated to principal autotopisms of long cycles. Application in Cryptography, Proceedings of Transgressive Computing 2006: a conference in honor of Jean Della Dora (2006), 213-230.

\bibitem{FalconAC} R. M. Falc\'on, Cycle structures of autotopisms of the Latin squares of order up to 11, {\em Ars Combinatoria} (in press). Available from http://arxiv.org/abs/0709.2973.

\bibitem{FalconHMAMS} R. M. Falc\'on, Designs based on the cycle structure of a Latin square autotopism, Proceedings of 1st Hispano-Moroccan Days on Applied Mathematics and Statistics (2008), 479-484.

\bibitem{FalconMartinJSC07} R. M. Falc\'on and J. Mart\'in-Morales, Gr\"obner bases and the number of Latin squares related to autotopisms of order $\leq$ 7, J. Symbolic Comput. 42 (2007), 1142-1154.

\bibitem{FalconMartinEACA08} R. M. Falc\'on and J. Mart\'\i n-Morales, The 3-dimensional planar assignment problem and the number of Latin squares related to an autotopism, Proceedings of XI Spanish Meeting on Computational Algebra and Applications (2008), 89-92.

\bibitem{Ghandehari05} M. Ghandehari, H. Hatami, E. S. Mahmoodian, On the size of the minimum critical set of a Latin square, Discrete Math. 293 (2005) 121-127.

\bibitem{Hardy18} G. H. Hardy, S. and Ramanujan, Asymptotic Formulae in Combinatory Analysis, Proc. London Math. Soc. 17 (1918), 75-115.

\bibitem{Hulpke11} A. Hulpke, P. Kaski and P. R. J. \"Osterg{\aa}rd, The number of Latin squares of order 11, Math. Comp. 80 (2011), 1197-1219.

\bibitem{Kerby10a} B. Kerby and J. D. H. Smith, Quasigroup automorphisms and symmetric group characters, Comment. Math. Univ. Carol. 51 (2010), 279-286.

\bibitem{Kerby10b} B. Kerby and J. D. H. Smith, Quasigroup automorphisms and the Norton-Stein complex, Proc. Amer. Math. Soc. 138, No.9 (2010), 3079-3088.

\bibitem{Kumar99} S. R. Kumar, A. Russell and R. Sundaram, Approximating Latin square extensions, Algorithmica 24 (1999), 128–138.

\bibitem{Laywine81} C. Laywine, An expression for the number of equivalence classes of Latin squares under row and column permutations, J. Combin. Theory Ser. A 30 (1981), 317-321.

\bibitem{Laywine85} C. Laywine and G. L. Mullen, Latin cubes and hypercubes of prime order, Fibonacci Quart. 23 (1985), 139-145.

\bibitem{McKay05}  B. D. McKay, A. Meynert and W. Myrvold. Small Latin Squares, Quasigroups and Loops, J. Comb. Designs, 15, No.2 (2007), 98-119.

\bibitem{McKay} B. D. McKay. http://cs.anu.edu.au/$\sim$bdm/data/latin.html.

\bibitem{Sade68} A. A. Sade, Autotopies des quasigroupes et des syst\`emes associatives, Arch. Math. 4, No. 1 (1968), 1-23.

\bibitem{Stones10} D. S. Stones, The parity of the number of quasigroups, Discrete Math. 310 (2010), 3033-3039.

\bibitem{Stones11} D. S. Stones, P. Vojt\v{e}chovsk\'y and I. Wanless, Cycle structure of autotopisms of quasigroups and Latin squares. Preprint avaliable from http://www.du.edu/media/documents/nsm/mathematics/preprints/ m1101.pdf

\bibitem{Wanless04} I. M. Wanless, Diagonally ciclic Latin squares, European J. Combin 25 (2004), 393-413.

\end{thebibliography}
\end{document}